\documentclass[12pt, twoside,a4paper]{article}
\pdfoutput=1
\usepackage[cp1251]{inputenc}
\usepackage{ifthen}
\usepackage{amsmath}
\usepackage{amsfonts}
\usepackage{latexsym}
\usepackage{amsthm}
\usepackage{amssymb}
\newcommand{\CopyName}{ V.\ M.\ Zhuravlov}
\newcommand{\NAME}{ V.\ M.\ Zhuravlov}
\newcommand{\Year}{2023}
\newcommand{\rightheadtext}{Predicates and terms from non-standard sequences}
     \newcounter{chapter}
     \newcounter{artpage}[chapter]
     \newcommand{\vs}{\vspace{.1in}}
     \newcommand{\vsk}{\vspace{.2in}}
\makeatletter
     \renewcommand{\@evenhead}{\footnotesize \ifthenelse{\value{artpage}=0}
     {\hfil}{\thepage\hfil \textsc {\leftmark} \hfil } }
     \renewcommand{\@oddhead}{\footnotesize\ifthenelse{\value{artpage}=0}
     {\hfil}{\hfil \textsc \rightmark \hfil \thepage} }
     \newcommand{\logo}{\baselineskip2pc \hbox to\hsize{\hfil\copyright\,\footnotesize
     \CopyName, \Year}}
     \renewcommand{\@oddfoot}{\ifthenelse{\value{artpage}=0}{\logo
     \refstepcounter{artpage}} {\hfil\refstepcounter{artpage}}}
     \renewcommand{\@evenfoot}{\ifthenelse{\value{artpage}=0}{\logo
     \refstepcounter{artpage}} {\hfil\refstepcounter{artpage}}}
     \renewcommand{\section}{\@startsection{section}{1}{0pt}{3.5ex plus
     1ex minus .2ex}{2.3ex plus 2.ex}{\large\hfil\textsc}}
     \vspace{3pt}\vs

\vs
\newcommand{\tit}{Predicates and terms from non-standard sequences}
\date{2023}
\usepackage{hyperref}
\usepackage{graphicx}
\graphicspath{{Images/}}
\begin{document}
\hfill
\vspace{0.3in}
\markboth{{\NAME}}{{\rightheadtext}}\begin{center} \textsc {\CopyName} \end{center}\begin{center} \renewcommand{\baselinestretch}{1.3}\bf {\tit} \end{center}
\vspace{20pt plus 0.5pt} {\abstract{\noindent
The article proposes a method for constructing non-standard theories based on terms from partially existing sequences of elements. The method is illustrated by the example of the theory of monoids.\newline
\textit{Predicates and terms from non-standard sequences, 2023, udc: 510.65 msc: 03H05\vspace{3pt}}\newline
\textit{Key words: monoid, term, ideal, category, isomorphizm, insignificant extension of the theory, non-standard theory}}
}\vsk
\tableofcontents
\section{Introduction}
What is a sequence of variables in its termal, logical sense? This is a certain subset of the Cartesian degree of the set of individuals for which a certain predicate can be fulfilled. Empirically (linguistically), this is a linear notation generated in a certain alphabet - which is what it is in the theory of Markov algorithms. And how can a sequence be defined algebraically without resorting to other constructions? It is striking that this product is in a certain free monoid, the unit of which is an empty symbol. One of the implementations of this approach are modern language constructs. But this definition is too narrow. It does not cover other interesting varieties of what might be called a sequence.\par
Monoids and semigroups are not only rather fundamental, but also rather poor constructions. Some improvement of these concepts is the concept of a category. The monoid is a trivial special case of the category. But perhaps the situation is not so trivial. In the article, we will pay attention to one interesting property of monoids and categories.\par
Let's try to proceed from the main "working" properties of sequences. One of the properties that make them such a versatile tool of thinking. In fact, sequences are what we call words. However, the meaningful word does not always exist.\newline

\smallskip
\section{Monoids and categories}
\textbf{Definition 1}: \textit{Let there be a monoid \textbf{M} with unit \textbf{E} and set \textbf{Q}, for which}:
$$\forall a(a\neq E)\forall b(b\neq E)\forall  c(c\neq E):(abc\in Q)\Rightarrow((ab\in Q)\vee(bc\in Q))$$
$$\forall x\forall y:(x\in Q)\Rightarrow((xy\in Q)\wedge(yx\in Q))$$
\textit{— then we will call any non-singular element from \textbf{(M-Q)} a chain (in this case, some chains can be equal, or consist of only one element-"link"). The set \textbf{Q} is called the big zero. If \textbf{Q} is singleton, then we call it associative zero and denote it as \textbf{0}.}\newline
Note: in quantifier expressions like: $\forall a:(a\neq E)...$ the statement is meant: "for all a not equal to \textit{\textbf{E…}}"\newline\par
Let me explain the meaning (the second of the conditions simply says that Q is a two-sided ideal). It is intuitively clear that if in a sequence of three elements one of the pairs does not form a chain, then the whole consecutive triplet will not be a chain either (well, the "hooks" are different and do not hook). Moreover, if an element is not a chain, then any of its possible "linking" with another element will also not be a chain. The result of a binary operation is unique for any pair of elements. An associative binary operation is unique for any sequence of elements (including single ones). A sequence (in its traditional sense) is formally characterized by uniqueness and associativity; there is one more, the main of its properties: any sequence exists if and only if any of its subsequences (any of its "continuous pieces") exist. It is easy to see that for our definition of a chain, all properties of sequences are satisfied by induction: the existence of a finite sequence is equivalent to the existence of all binary sequences included in it (according to Markov). But ... a simple feature - so defined (chain) sequences do not always exist. But the existence of a certain finite chain implies the existence of all its subchains (in the order of writing). And if we abstract from the fourth fundamental property of sequences, from the existence of a sequence of any chains of symbols, then we get our definition of a chain.\par
So, belonging to \textit{\textbf{Q}} means "not being a sequence". A chain is a non-standard definition of a sequence. \textit{\textbf{Q}} is a collection of "non-strings", non-existent sequences. And what \textit{\textbf{Q}} in the monoid? This is a two-sided associative (as one would like to call it) ideal (see the first property in \textbf{Definition 1}). Which is an approximation to some prime ideal \textit{\textbf{P}} (in the it: $\forall x\forall y: (xy\in P)\Leftrightarrow[(x\in P)\vee(y\in P)]$, where the right implication is the definition of a completely divisible set). Since the first property in \textbf{Definition 1} is a certain weakening of the simplicity condition (divisibility of the product of triples of factors). The true lemma is:\newline\par
\textbf{Lemma 1}. \textit{Every prime ideal is associative. The reverse is not true. The definition of an associative ideal is consistent.}\newline\par
Proof: $\forall a\forall b\forall c:[(ab\in P)\vee(bc\in P)]\Leftrightarrow(abc\in P)$,— by the definition of the simplicity of an ideal. Let us give an example of an associative, but not a simple, ideal. Consider a monoid in which there are no commuting products without the participation of unity or the power of the element. And let this monoid have a set of generators. Then all elements that are multiples of some generatrix constitute a simple ideal, not necessarily an associative one. However, elements that are multiples of the product of two generators will constitute an associative, but not a simple, ideal (the product of this very pair of generators is indivisible within this set). However, we will see other examples later. Since associative ideals exist, their definition is consistent. The proof is over.\newline\par
And what will happen when the monoid \textit{\textbf{M}} is factorized by the identity relation for its associative ideal: \textit{\textbf{(Q=0)}} (we mean the identification of all elements of the ideal, which results in one of its elements, denoted as \textbf{0})? We get a monoid \textit{\textbf{M0}} with an associative zero:
$$(abc=0)\Rightarrow[(ab=0)\vee(bc=0)]$$
$$x0=0x=0$$
— onto which M is mapped epimorphically, but the entire structure of non-zero elements remains unchanged. At the same time, the identification: \textit{\textbf{(Q = E)}} will lead to the triviality of the entire monoid: \textit{\textbf{(M = E)}}, since \textit{\textbf{Q}} is an ideal. More precisely, it is true:\newline\par
\textbf{Lemma 2:}  \textit{Co-decart square (it is commutative and the pair of arrows in the upper right corner is an amalgam of the pair of arrows from the lower left):\begin{center}\includegraphics{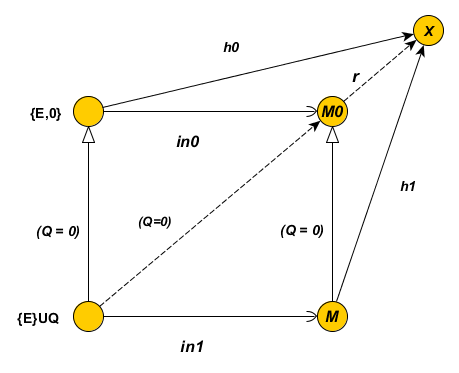}\newline\end{center}
where:}\par
\textit{\textbf{E} — is the null object of the category of monoids (isomorphic to the unity of all monoids);}\par
$(\{E\}\cup Q)$ \textit{is the monoid consisting of the ideal \textbf{Q} and the identity \textbf{E};}\par
$(\{E,0\})$ \textit{is a monoid consisting of \textbf{0} and \textbf{1} (similarly);}\par
\textit{\textbf{(E=Q), (E=0)}, and \textbf{(Q=0)} — these relations here represent the corresponding epimorphisms;}\newline\par
The proof is clear enough. Let we have a pair of monoid homomorphisms $h1:M\rightarrow X$ and $h0:\{E,0\}\rightarrow X$, such that
$$in\ast h1=\overrightarrow{(Q=0)}\ast h0$$
( $\overrightarrow{(Q=0)}$ means epimorphism $(\{E\}\cup Q)\ast (\{E,0\})$, induced by the ratio: $(Q=0)$ ). This commutativity means that \textit{\textbf{h0}} and \textit{\textbf{h1}} transform \textbf{0} and \textit{\textbf{Q}}  into the same element: $h0=h1$; define the mapping $r:M0\rightarrow X$,
$$r(0)=h0(0)$$
$$(x\not\in Q)\Rightarrow(r(x)=h1(x))$$
Then the image of \textit{\textbf{Q}} will be the only idempotent element belonging to the monoid \textit{\textbf{X}}, which fixes the only limit homomorphism $r:M0\rightarrow X$, (under which \textbf{0}  passes into the above idempotent), commutatively "passing through itself" the original pair of arrows \textit{\textbf{(h1, h0)}}. We get the complete commutativity of the diagram (including internal triangles). The proof is over.\newline\par
Now that we have a little understanding of the position of associative ideals in monoid theory, let's move on to studying the relationship they have to category theory.\par
We construct a "non-standard" model of a monoid with an everywhere defined product, but with the properties of a partial product. Moreover, the belonging of an element to an associative ideal will model its "non-existence". Indeed, the defining property of the associativity of an ideal can be rephrased as follows:
$$\forall a(A\neq E)\forall b(b\neq E)\forall c(c\neq E):((ab\not\in Q)\wedge(bc\not\in Q))\Rightarrow(abc\not\in Q)$$
Those — if the products to the left of the implication sign exist, then the product to the right also exists. And this means that the law of composition is associative not only in general on the entire monoid, but also where it "exists" relative to the zero ideal (consisting of elements, as it were, "non-existent"). The associativity of the law of composition is necessary (first of all) precisely for the uniqueness of the product of several elements (brackets inside the sequence of factors can be placed arbitrarily, despite the binary nature of the original operation). And she, this associativity, is one of the axioms of the theory of monoids. There is another axiom that asserts the existence of a unit. But we want to model the very concept of existence.\par
\textit{Therefore, in any sentence, the mention of the product \textbf{(xy)} must be preceded by the implication: $(xy\not\in Q)\Rightarrow ...$} (which must be introduced into the axiom of the existence of a unit).\par
\textit{The product, which does not always exist, becomes local with respect to each multiplicand} (and such locality must be consistent "right" and "left"). \textit{The existence of each individual element is determined by its existence as a multiplicand. But the multiplicand can be not one element, but a set of elements that have a given multiplicative property.} This set must be localized with respect to other elements, the product with which is involved in the formulation of this property. Moreover, this set can also be empty. For example, a unit whose existence axiom is preceded by the above implication may in fact not exist at all. \textit{Therefore, in any sentence, if you postulate the existence of \textbf{x}, such that ... and then the product \textbf{(xy)} follows, then your axiom must be preceded by a quantifier prefix:} $\forall x\exists y:...(x\ast y)...$, \textit{where} $(x\ast y)$ \textit{is a \textbf{non-empty} term, actually \textbf{defined} for the pair \textbf{(x,y)}.}\par
Consequently, the axiom of the existence of a unit must undergo the following procedure, which we call the localization of the unit of a monoid with respect to the associative ideal \textit{\textbf{Q}}. And, as we have already said, we will consider non-belonging to an associative ideal to be a synonym for existence. Then we have:\newline\par
\textbf{Definition 2:} \textit{We call the elements \textbf{e} of our monoid small, or local units, if:}
$$L(e)\equiv\forall x:[(ex\not\in Q)\Rightarrow(ex=x)]\wedge[(xe\not\in Q)\Rightarrow(xe=x)]$$
Let's add an (insignificant!) couple of axioms:
$$\forall x(x\not\in Q)\exists e:L(e)\wedge(ex\not\in Q)$$
$$\forall x(x\not\in Q)\exists e:L(e)\wedge(xe\not\in Q)$$
\textit{(in exactly the same way we will designate and call \textbf{e} and \textbf{L(e)}) unit arrows in the category, i.e. its objects (after all, they have the same properties); allowing "liberty of speech", we will also use the notation:} "$e\in L$", \textit{meaning belonging to a set of local units.}\newline\par
In what follows, by \textit{\textbf{M1}} we mean the monoid M extended by the collection of its small units, i.e. $M1\equiv M\cup L$. There is a monomorphic embedding: $M\hookrightarrow M1$. In the metamathematical sense, the extension of M to \textit{\textbf{M1}} is inessential (just as the extension of a semigroup to a monoid is also inessential). We simply introduce some limiting elements, regardless of whether they existed before. It is clear that: \textit{\textbf{M01=M10}}. It is also clear that small units do not have the uniqueness of a "large" unit of a monoid.\newline\par
\textbf{Theorem 1:} \textit{For every pair \textbf{(M, Q)} there is a unique category that maps epi-mono-morphically onto the set} $M1-(\{E\}\cup Q)$.\newline\par
Proof. This category is the set itself. We define a partially defined operation: $a\ast b\equiv a\cdot b|\not\in Q$, where "dot" means the product in the monoid, and "asterisk" is this product where it "exists", i.e. — is not included in \textit{\textbf{Q}}. We have an identical mapping (one-to-one by definition): $id:(M1-(\{E\}\cup Q))\rightarrow(M1-(\{E\}\cup Q))$, with the homomorphism property: $id(a\ast b)=a\cdot b$. In this case, small units are objects of the category (it is known that the unit arrows of a category are its objects, that is, what determines \textit{\textbf{cod}} and \textit{\textbf{dom}} each arrow; in fact, they are the limits and colimits of a diagram consisting of one arrow). We get an epi-monomorphism, but not an isomorphism — the properties of the compositions are still so different that "isomorphism" works only in one direction. The proof is over.\newline\par
It is noteworthy that the triple appears here: $(E, M1-(\{E\}\cup Q), 0)$ (where \textbf{0} is the big zero of the monoid), whose categorical analogue we will now define.
\section{Categories and monoids}
\textbf{Definition 3:} \textit{For each category, we define a triple $K1\equiv(\widehat{E},K,\widehat{\emptyset})$ with all intra-category operations and a composition of \textbf{K} elements extended with additional properties:}
$$\forall a(a\in K):(\widehat{E}\ast a=a\ast\widehat{E}=a)\wedge(\widehat{\emptyset}\ast a=a\ast\widehat{\emptyset}=\widehat{\emptyset})$$
$$\widehat{E}\ast\widehat{\emptyset}=\widehat{\emptyset}\ast\widehat{E}=\widehat{\emptyset}$$
$$\forall a(a\in K)\forall b(b\in K):(\exists c(\in K):c=a\cdot b)\Longleftrightarrow (a\cdot b=c)$$
$$\forall a(a\in K)\forall b(b\in K):(\nexists c(\in K):c=a\cdot b)\Longleftrightarrow (a\cdot b=\widehat{\emptyset})$$
\textit{where: \textbf{a, b} are any arrows from \textbf{K}, including single ones, i.e. its objects (idempotent arrows like: $\circlearrowleft$); "dot" denotes within-category composition, and "asterisk" denotes the extended composition that we define.}\newline\par
If you like, you can consider $\widehat{E}$ the identity functor for \textit{\textbf{K}} (or its "big" unit) , and $\widehat{\emptyset}$ its erasing functor (or the empty set). At the same time, it is obvious that we will get an insignificant extension of the theory \textit{\textbf{K}}: any sentence, limited by its elements, retains its truth or falsity.\newline\par
\textbf{Theorem 2:} \textit{Every category \textbf{K} embeds monomorphically into some monoid (not necessarily strict) \textbf{M} with an associative ideal \textbf{Q}. Any true or false sentence of the theory \textbf{K1}, limited only by elements from \textbf{K}, corresponds to the same (true or false) sentence of the theory \textbf{K}. And vice versa.}\newline\par
Proof. We define a mapping: $\mu: K\rightarrow K1$, such that
$$\forall a(a\in K):\mu(a)=a$$
$$\forall a(a\in K)\forall b(b\in K):[(\exists c(c\in K):c=a\cdot b)\Rightarrow\mu(a\cdot b=a\ast b)]\wedge$$
$$\wedge[(\nexists c(c\in K):c=a\cdot b)\Rightarrow\mu(a\cdot b)=\widehat{\emptyset}]$$
(in the first formula, "liberty of speech" is allowed: the same a is an element in a different hierarchy of sets — the essence of the question does not change from this)\par
The construction \textit{\textbf{K1}} constructed in \textbf{Definition 3} is a monoid with a large zero: there is a unit, the binary composition is everywhere associative, and $\widehat{\emptyset}$, this is the associative ideal of one element. The map $\mu$ bounded by the images of elements from \textit{\textbf{K}} (identical to these elements) is the identical isofunctor (in fact, it can be identified with $\widehat{E}$, as we have already noted). The mapping $\mu$ defines a one-to-one correspondence of elements that exist in \textit{\textbf{K}}. A note about associativity. According to the categorical axiom:
$$(\exists z1:z1=a\cdot b)\wedge(\exists z2:z2=b\cdot c)\Longrightarrow\exists z:z=abc$$
But this is equivalent to the associativity of the ideal $\widehat{\emptyset}$ (and it is an ideal, since it is an idempotent multiplicative zero):
$$(a\ast b\ast c=\widehat{\emptyset})\Longrightarrow[(a\ast b=\widehat{\emptyset})\vee(b\ast c=\widehat{\emptyset})]$$
We have already spoken about the inessentiality of the extension \textit{\textbf{K}}. It is a consequence of the monomorphism of the embedding of \textit{\textbf{K}} in \textit{\textbf{K1}}. Therefore the theorem is proved.\newline\par
And it is clear that the procedure used by us in the proofs of both theorems is completely reversible: weakening the monoid obtained from the category by rejecting the presence of small ones and replacing the large zero with the associative ideal \textit{\textbf{Q}}, we arrive at the monoid with which we started.\par
In fact, we have also proved the following theorem:\newline\par
\textbf{Theorem 3:} \textit{Every category is a complete model of some monoid with an associative zero. Conversely, any such monoid completely models some category. Category theory is logically equivalent to the theory of monoids with an associative zero. Both are non-essential extensions of each other.}\newline\par
The proof is quite simple. Let us repeat his main theses. Objects of an arbitrary category are identified with single arrows. And then we will look at the category as a set of arrows with a partially defined composition. Let's make the composition completely defined by adding a certain object \textbf{0} to the category, to which we will equate all non-existent compositions of arrows. We will also introduce the unit functor of the category into our construction as \textit{\textbf{E}}, the identity of the monoid. Since categorical composition is associative where it exists, we indeed include our category isomorphically in a minimal monoid with a large zero. Or, if you like, into the semigroup of such a monoid.\par
And vice versa, non-trivial (not \textit{\textbf{E}} and not \textbf{0}) elements of any monoid with an associative zero can be identified with arrows of a certain category by declaring zero products non-existent. And by associating E with its identity functor. This will also lead to an isomorphic immersion of the monoid into some minimal category.\par
As for logical equivalence, such artifacts as \textit{\textbf{E}} and category objects (unlike its arrows) do not carry new information and are only insignificant extensions of the corresponding theories (and in category theory, \textbf{0} will not be a significant extension). The theorem has been proven.\newline\par
Above, we linked the objects (single arrows) of the category and the elements of the monoid. Actually, the unit of the monoid becomes a set of small units when it "transforms" into a category. And the set of objects of the category becomes the unit of the monoid in the reverse process. If we identify the elements with the operations that they perform (due to multiplication) on other elements, then the "big" unit \textit{\textbf{E}} takes the element \textit{\textbf{x}} from \textit{\textbf{M}} into itself when acting from the left and right. When mapped onto some small unit (i.e., onto a single arrow, object \textit{\textbf{K}}), it will translate into itself a certain arrow (acting from the left or right). This arrow will reappear at \textit{\textbf{x}}. We get a commutative diagram. At the same time, such artifacts as: $(0,Q,\widehat{0},...)$ — in the case of a category are sets of pairs of arrows with a non-existent composition, and in the case of a monoid, they are sets of "non-existent" products.\par
Let us describe the situation even more precisely—there is a homomorphic relation R between an arbitrary category and some monoid with an associative ideal. Those,— everything is the same as with the usual homomorphism, but without the condition of functional uniqueness:
$$R\subseteq (M1\times K1)$$
$$\forall a(a\in M)\exists x(x\in K1):R(a,x)$$
$$(R(a,x)\wedge R(b,x))\Rightarrow R(ab,xy)$$
Now let's define the properties:
$$R(E,\widehat{E})\wedge[\exists!x(x\in K1):R(E,x)]$$
$$R(Q,o)\wedge[\exists!x(x\in K1):R(Q,x)]$$
— the unit of the monoid always corresponds to the "large unit" of the category; the associative ideal is always completely transformed into the element \textit{\textbf{o}}.
$$\forall a(a\in (M1-\{E\}\cup Q))\exists!x(x\in K):R(a,x)$$
$$\forall x(x\in K)\exists!a(a\in (M1-\{E\}\cup Q)):R(a,x)$$
— an extended monoid without unity and zero ideal, is in one-to-one correspondence with category arrows;
$$\forall e(e\in M1):L(e)\Rightarrow L(R(L(e)))$$
$$\forall e(e\in K):L(e)\Rightarrow L(R(L(e)))$$
— small units of the monoid and the category, in turn, also one-to-one correspond to each other.\newline
The structure of relation \textit{\textbf{R}} is as follows:\begin{center}\includegraphics{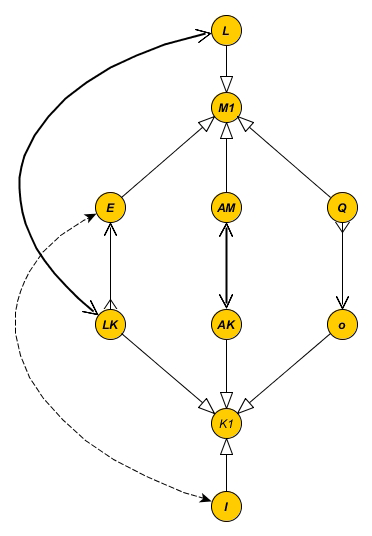}\newline\end{center}
— where two-sided arrows denote isomorphisms; arrows with transparent "tips" — inclusions in the direct sum (combinations of collections); arrows with a branched base are epimorphic mappings into a single element.\par
And these correspondences preserve the operation of composition where they exist. Moreover, if we also introduce a set of non-existent categorical compositions into consideration, we get something isomorphic to the zero ideal. Based on this ideal, it is possible to create certain structures that, in a certain sense, are additional to categorical ones. The category theory itself can be immersed in the category of monoids.\par
This is the actual identification of categories with some type of monoids. And not always obvious immersion of the theory of monoids in the theory of categories. I think that such a representation can be more productive than the trivial identification of a monoid with a category of one object.
\section{Partially existing sequences}
Earlier we talked about the fundamental properties of sequences:\par
•	uniqueness (i.e., the associativity of the operation of joining sequences); the existence of those and only those sequences that consist only of existing subsequences (possibility of arbitrary reduction on the left and right).\par
•	for any existing sequence, there is a pair of elements that, in combination with it from the left (right), will create existing sequences (the possibility of infinite extension of the sequence to the right and left).\par
•	the presence of an empty sequence (the result of assigning to which does not change anything). In a monoid, the role of the empty sequence is played by the unit. In the category — single arrows (objects).\par
•	the existence of an element is equivalent to the existence of its sequence with unit length.\newline\par
From this it follows (by induction on any finite length) that the logical predicate of the existence of sequences is completely determined by its binary part.\par
The theory should be modeled in non-standard sequences only on existing elements.\par
In the monoid \textit{\textbf{M}}, we consider as existing any element that does not belong to the given associative ideal. Then we obtain a model for the existence of sequences — all four properties are satisfied: the only difference is that not all sequences exist. Indeed, let's define the predicate: $ex(x)\equiv(x\not\in Q)$. Then, according to \textbf{definition 1}:
$$(ex(ab)\wedge ex(bc))\Leftrightarrow ex(abc)$$
the predicate \textit{\textbf{ex}} for any existing string is satisfied if and only if it is satisfied for all its binary substrings, which is proved by induction; there are chains of unit length, because the ideal does not include all the elements of the monoid; moreover, \textit{\textbf{E}} should be identified with the empty string.\par
Therefore, an \textit{\textbf{ex}}-predicate defined on a free monoid has all the properties of a partial existence quantifier for sequences. And in a free monoid, products are one-to-one with sequences of elements. At the same time, all relations between elements of any monoid correspond one-to-one to relations between category arrows.\par
There is one more remark to which it is possible to pay or not to pay attention. According to the fourth property of sequences, the associative ideal of a monoid must not be totally divisible by an element that does not belong to it:
$$(x\not\in Q)\Longrightarrow\exists y:(y\not\in Q)\wedge(xy\not\in Q)$$
$$(x\not\in Q)\Longrightarrow\exists y:(y\not\in Q)\wedge(yx\not\in Q)$$
 — this is a very weak analogue of the simplicity of the ideal. This analogue also makes the ideal "a little more"...\par
Let us pay attention to the fact that the existence of a composition of arrows (or elements of a monoid), being partial, is localized with respect to factors. And localization is performed at a more fundamental level — on sequences of elements. This is obvious enough.
Let us turn to a more general theme, which all the previous ones suggest. Let's define a procedure:\newline\par
\textbf{Definition 4}: \textit{On all finite sequences, we define the existence predicate:}
$$\forall a\forall b\forall c:ex(a,b)\wedge ex(b,c)\Leftrightarrow ex(a,b,c)$$
$$\exists a\exists b:ex(a,b)$$
$$\exists a: ex(a)$$
$$\forall x: ex(x)\Rightarrow\exists y\exists z: ex(y)\wedge ex(z)\wedge ex(x,y)\wedge ex(z,x)$$\newline
Such axioms say that this is actually a binary predicate, obviously extendable to all finite sequences (similar to the associative product operation) and non-empty on all finite sequences. Such a predicate can be any binary predicate, which for any element is "non-empty" on the left and on the right. Sometimes, as in the case of a monoid with an associative ideal, such a predicate can already be naturally defined in theory. \textit{And then, in all the axioms of the theory, we make the substitution:}
$$\rho(x_{1},...,x_{n})\longmapsto\bigwedge_{j\in (1,...,n-1)}ex(x_{1},x_{j+1})\Rightarrow\rho(x_{1},...,x_{n})$$
\textit{where $\rho$ is an arbitrary predicate from the axioms of the theory. In addition, predicates with constants are replaced (one predicate for two):}
$$\rho(...,x,c,...)\longmapsto\rho(...,x,\tau 1(x),...)$$
$$\rho(...,c,y,...)\longmapsto\rho(...,\tau 2(y),y,...)$$
\textit{where $\tau 1$ and $\tau 2$ are one-place terms replacing zero-place ones (that is, they are constants of the theory; it is clear that such a term implies the existence of the corresponding element).}\par
It is this procedure (but in somewhat different terms) that we carried out in the theory of a monoid with an associative ideal, and arrived at category theory. This indicates the correctness of such substitutions.\par
It should be noted that the "binary non-emptiness" of the ex-predicate (the last of its properties) also leads to the "non-emptiness of the implication" in the first replacement: i.e., a sequence with this property exists in the standard sense:$\exists(x_{1},...,x_{n}):\rho(x_{1},...,x_{n})$.\par
In addition, the functional nature of the "right" and "left" terms from the second and third substitutions localizes the constants with respect to each element (which does not at all mean that they are unique for this element).\par
\textit{And the following suggestion seems very likely:}\newline\par
\textbf{Hypothesis.} \textit{Any mathematical theory can be constructed in a logic in which terms and sequences are localized with respect to individuals according to the procedure from \textbf{Definition 4}. The existence of sequences and the results of operations on them, in this case, is modeled non-standard.}\newline
The argument preceding the proof of \textbf{Theorem 1} also speaks in favor of the conjecture.\newline\par
We have seen that the very existence of any sequence of elements can be localized with respect to each of its elements. And this leads to a transformation of any mathematical theory that uses terms, to a change in functions, operations, and even relations. The radical nature of such changes is clearly seen in the example of monoids and categories.\par
\textit{\textbf{The category is a locally thermal monoid.}} Locally thermal theories are non-standard theories, with non-standard sequences and non-standard existence. Which is very similar to non-standard according to Robinson. All this uses certain modifications of ultrafilters and compactness theorems at various metatheoretical levels. And it allows you to somehow simulate contradictions and approach the solution of paradoxes.\newline
\begin {center}\textbf{Bibliography}\end{center}

\noindent [1] Abraham Rodinson,University California,\par
\textit{Introduction to model theory fnd to the metamathematics of algebra} (1963).\newline

\noindent [1] R. Goldblatt. North-Holland Publishing Company: Amsterdam, New York, Oxford.\par
\textit{Topoi The categorial analysis of logic}, (1979).\newline

\noindent [2] Haskell B. Curry, \textit{Foundations of Manthematical logic},\par
McGraw-HILL Book Company, INC. (1984).\newline

\noindent [3] Helena Rasiowa and Roman Sikorski, \textit{The Mathemayics of Metamathematics},\par
Panstwowe Wydawnlctwo Naukowe Warszawa, (1963).\newline

\noindent [4] V. A. Uspensky. \textit{What is non-standard analysis?}, Nauka Publishing House, (1987).

\end{document}